\documentclass[11pt]{amsart}
\title{Products and selection principles.}
\author{Liljana Babinkostova and Marion Scheepers}
\date{}
\usepackage{amssymb}
\newcommand{\sone}{{\sf S}_1}
\newcommand{\gone}{{\sf G}_1}
\newcommand{\sfin}{{\sf S}_{fin}}
\newcommand{\gfin}{{\sf G}_{fin}}
\newcommand{\Sc}{{\sf S}_c}

\newcommand{\naturals}{{\mathbb N}}
\newcommand{\reals}{{\mathbb R}}

\newtheorem{theorem}{Theorem}
\newtheorem{lemma}[theorem]{Lemma}
\newtheorem{proposition}[theorem]{Proposition}	
\newtheorem{corollary}[theorem]{Corollary}	
\newtheorem{problem}{Problem}
\newtheorem{conjecture}{Conjecture}

\begin{document}
\maketitle
\begin{abstract} 
We study when the product of separable metric spaces has the selective screenability property, the Menger property, or the Rothberger property. Our results imply The product of a Lusin set and 
\begin{enumerate}
\item{a Sierpinski set always has the Menger property (Corollary \ref{LusinSierpinski});}
\item{a $\gamma$-set always has the Rothberger property (Corollary \ref{gammaLusin}).}
\end{enumerate}
\end{abstract}

All topological spaces considered in this paper are assumed to be separable metric spaces. Some cited theorems apply to more general spaces, as the reader could verify by consulting the appropriate references. 

Let $\mathcal{A}$ and $\mathcal{B}$ be given families of collections of subsets of some set $S$. The following selection principle was introduced in \cite{lbthesis}:
\begin{quote}
$\Sc(\mathcal{A},\mathcal{B})$: For each sequence $(O_m:m<\infty)$ of elements of $\mathcal{A}$ there is a sequence $(T_m:m<\infty)$ with each $T_m$ a pairwise disjoint family refining $O_m$, and $\cup\{T_m:m<\infty\}\in\mathcal{B}$.
\end{quote}

The special case of $\Sc(\mathcal{A},\mathcal{B})$ when $\mathcal{A}=\mathcal{B}=\mathcal{O}$, the collection of open covers of a topological space, was introduced in \cite{ag} by Addis and Gresham. This property is related to the theory of covering dimension. One of the interesting questions about it, due to D. Rohm, asks when the product of two spaces with $\Sc(\mathcal{O},\mathcal{O})$ again has this property. The best known result regarding this question is due to \cite{hy} and \cite{rohm}:
\begin{theorem}[Hattori-Yamada, Rohm]\label{scproduct} Let $X$ and $Y$ be topological spaces satisfying $\Sc(\mathcal{O},\mathcal{O})$. If $X$ is $\sigma$-compact, then $X\times Y$ has the property $\Sc(\mathcal{O},\mathcal{O})$.
\end{theorem}
The special case when $\mathcal{A} = \mathcal{T}$, the collection of two-element open covers of a space, and $\mathcal{B} = \mathcal{O}$, was introduced by Aleksandroff and is known as \emph{weak infinite dimensionality}. This is not the original definition of weak infinite dimensionality, but Rohm has shown that $\Sc(\mathcal{T},\mathcal{O})$ is equivalent to weak infinite dimensionality. It is unknown whether $\Sc(\mathcal{T},\mathcal{O})$ is equivalent to $\Sc(\mathcal{O},\mathcal{O})$. Various alternatives of the hypothesis that one of the spaces be $\sigma$-compact have been investigated in attempts to generalize Theorem \ref{scproduct}. To explain these we now recall two more selection principles from \cite{coc1}: 

\begin{quote}
$\sfin(\mathcal{A},\mathcal{B})$: For each sequence $(O_m:m<\infty)$ of elements of $\mathcal{A}$ there is a sequence $(T_m:m<\infty)$ with each $T_m$ a finite subset of $O_m$, and $\cup\{T_m:m<\infty\}\in\mathcal{B}$.
\end{quote}
The \emph{Menger property} is $\sfin(\mathcal{A},\mathcal{B})$ when $\mathcal{A}=\mathcal{B}=\mathcal{O}$. Hurewicz introduced it in \cite{wh}. 

\begin{quote}
$\sone(\mathcal{A},\mathcal{B})$: For each sequence $(O_m:m<\infty)$ of elements of $\mathcal{A}$ there is a sequence $(T_m:m<\infty)$ with each $T_m\in O_m$, and $\{T_m:m<\infty\}\in\mathcal{B}$.
\end{quote}
$\sone(\mathcal{O},\mathcal{O})$, introduced in \cite{rothb} by Rothberger, is known as the \emph{Rothberger property}. 

The following two results are probably the best known ones regarding products having $\sfin(\mathcal{O},\mathcal{O})$ or $\sone(\mathcal{O},\mathcal{O})$:
\begin{theorem}[Folklore]\label{sigmacompact} If $X$ is $\sigma$-compact and $Y$ has $\sfin(\mathcal{O},\mathcal{O})$, then $X\times Y$ has $\sfin(\mathcal{O},\mathcal{O})$.
\end{theorem}

\begin{theorem}[Folklore]\label{countable} If $X$ is countable space and $Y$ has $\sone(\mathcal{O},\mathcal{O})$, then $X\times Y$ has $\sone(\mathcal{O},\mathcal{O})$. 
\end{theorem}

We will investigate Rohm's question for $\Sc(\mathcal{O},\mathcal{O})$ as well as its analogues for $\sfin(\mathcal{O},\mathcal{O})$ and $\sone(\mathcal{O},\mathcal{O})$. In Section 1 we briefly survey some known limitations. In Section 2 we point out a common thread to Theorems \ref{scproduct}, \ref{sigmacompact} and \ref{countable}. In Section 3 we prove a few new results. A consequence of one of our results is that the product of a Sierpinski set and a Lusin set has Menger's property. In Section 4 we state a conjecture and an open problem about possible generalizations of Theorem \ref{pot}. 

\section{Limitations on the factors in products.} 

Various types of open covers are relevant to this discussion. Here are some of them: An open cover $\mathcal{U}$ of a topological space is:
\begin{itemize}
\item{A \emph{large} cover if for each $x\in X$ the set$\{U\in\mathcal{U}:x\in U\}$ is infinite.}
\item{An $\omega$-cover if $X\not\in \mathcal{U}$, and for each finite set $F\subset X$ there is a $U\in\mathcal{U}$ with $F\subset U$.}
\item{\emph{groupable} if there is a disjoint partition $\mathcal{U} = \cup_{n\in\naturals}\mathcal{U}_n$ into finite sets $\mathcal{U}_n$ such that each element of $X$ is in all but finitely many of the sets $\cup\mathcal{U}_n$.}
\item{A $\gamma$ cover if it is an infinite cover, and each infinite subset of it is a cover of the space.}
\end{itemize}
The symbol $\Lambda$ denotes the collection of large covers, $\Omega$ denotes the collection of $\omega$-covers, $\mathcal{O}^{gp}$ denotes the collection of groupable open covers and $\Gamma$ denotes the collection of $\gamma$ covers of $X$. Note that $\mathcal{O}^{gp} = \Lambda^{gp}$. In \cite{coc7} it was shown that another property introduced by Hurewicz in \cite{wh} is equivalent to $\sfin(\Omega,\mathcal{O}^{gp})$. This property is known as the Hurewicz property.

Regarding Theorem \ref{scproduct}: R. Pol showed in \cite{rpol} that weakening $\sigma$-compactness to $\sfin(\mathcal{O},\mathcal{O})$ does not generalize this theorem, even if the product of the factor spaces is assumed to have the Menger property, $\sfin(\mathcal{O},\mathcal{O})$:
\begin{theorem}[R. Pol]\label{rpoltheorem} Assume the Continuum Hypothesis. Then there is for each positive integer $n$ a separable metric space $X$ such that: 
\begin{itemize}
\item{$X^n$ has both $\sfin(\mathcal{O},\mathcal{O})$ and $\Sc(\mathcal{O},\mathcal{O})$, and}
\item{$X^{n+1}$ has $\sfin(\mathcal{O},\mathcal{O})$ but not $\Sc(\mathcal{O},\mathcal{O})$.}
\end{itemize}
\end{theorem}

It is shown in \cite{lbsc} that if the Menger property is strengthened to the Hurewicz property (introduced in \cite{wh} by Hurewicz and still weaker than $\sigma$-compactness), then a positive result is obtained for finite powers:
\begin{theorem}[\cite{lbsc}]\label{babinkostovath} 
If $X$ is a separable metric space such that $X^n$ has $\sfin(\Omega,\mathcal{O}^{gp})$ and $X$ has $\Sc(\mathcal{O},\mathcal{O})$, then $X^{n}$ has $\Sc(\mathcal{O},\mathcal{O})$.
\end{theorem} 

E. Pol investigated another alternative: Replace $\sigma$-compactness of one of the factor spaces with zerodimensionality, a strengthening of $\Sc(\mathcal{O},\mathcal{O})$. E. Pol showed in \cite{epol86} that the product theorem fails also for this alternative:
\begin{theorem}[E. Pol]\label{epolth} There exist separable metric spaces $X$ and $Y$ such that $X$ is zerodimensional and $Y$ has property $\Sc(\mathcal{O},\mathcal{O})$ but $X\times Y$ does not have property $\Sc(\mathcal{O},\mathcal{O})$. If the Continuum Hypothesis is assumed, the space $X$ can be taken to be the set of irrational numbers.
\end{theorem}

Motivated by this line of inquiry we drop $\sigma$-compactness of one factor, but strengthen $\Sc(\mathcal{O},\mathcal{O})$ in both factors. A positive result is obtained in Theorem \ref{pot} (1). 

Regarding Theorems \ref{sigmacompact} and \ref{countable}: 
A space which has the selection property $\sone(\Omega,\Gamma)$ also has the property $\sone(\mathcal{O},\mathcal{O})$, and thus the property $\sfin(\mathcal{O},\mathcal{O})$ - see \cite{coc2} and \cite{coc1}. Gerlits and Nagy introduced the property $\sone(\Omega,\Gamma)$ in \cite{GN}. Sets of real numbers with this property are also called $\gamma$-sets. It is consistent relative to the consistency of {\sf ZFC} that all separable metric spaces with $\sone(\Omega,\Gamma)$ are countable. In \cite{GM} Galvin and Miller showed that under appropriate hypotheses there are uncountable sets of real numbers which have $\sone(\Omega,\Gamma)$. Using their techniques one can show that the Continuum Hypothesis implies that there are sets $X$ and $Y$ of real numbers, each with $\sone(\Omega,\Gamma)$, for which $X\times Y$ does not have $\sfin(\mathcal{O},\mathcal{O})$. In \cite{todorcevic} Todor\v{c}evic gave {\sf ZFC} examples of (nonmetrizable) spaces $X$ and $Y$ which have $\sone(\Omega,\Gamma)$, but $X\times Y$ does not have $\sfin(\mathcal{O},\mathcal{O})$.

Again we drop the $\sigma$-compactness (or countability) of one factor in Theorems \ref{sigmacompact} and \ref{countable} but strengthen the corresponding selection principles in both factors. Positive results are obtained in Theorem \ref{pot} (2) and (3). 

\section{The game-theoretic connection.}

Selection principles have natural games associated with them. For an ordinal number $\alpha$, define:
\begin{quote}
$\gfin^{\alpha}(\mathcal{A},\mathcal{B})$: This is a game with $\alpha$ innings. In inning $\gamma<\alpha$ ONE first chooses a set $O_{\gamma}\in\mathcal{A}$. Then TWO responds with a finite set $T_{\gamma}\subseteq O_{\gamma}$. TWO wins a play
$O_0,\, T_0,\, \cdots,\, O_{\gamma},\, T_{\gamma},\, \cdots$ if $\bigcup\{T_{\gamma}:\gamma <\alpha\}\in\mathcal{B}$. Else, ONE wins.
\end{quote}
And similarly define:
\begin{quote}
$\gone^{\alpha}(\mathcal{A},\mathcal{B})$: This is a game with $\alpha$ innings. In inning $\gamma<\alpha$, ONE first chooses a set $O_{\gamma}\in\mathcal{A}$. Then TWO responds with a $T_{\gamma}\in O_{\gamma}$. TWO wins a play
$
  O_0,\, T_0,\, \cdots,\, O_{\gamma},\, T_{\gamma},\, \cdots
$
if $\{T_{\gamma}:\gamma <\alpha\}\in\mathcal{B}$. Else, ONE wins.
\end{quote}
It is clear that if for some countable ordinal $\alpha$ ONE has no winning strategy in the game $\gfin^{\alpha}(\mathcal{A},\mathcal{B})$, then $\sfin(\mathcal{A},\mathcal{B})$ holds. The converse of this is not always true. The same remarks apply to the game $\gone^{\alpha}(\mathcal{A},\mathcal{B})$. Also note that if TWO has a winning strategy in $\gone^{\alpha}(\mathcal{A},\mathcal{B})$, then TWO has a winning strategy in $\gfin^{\alpha}(\mathcal{A},\mathcal{B})$.

Existence of a winning strategy for TWO imposes structure on the underlying space. The two best known classical results in this connection are:
\begin{theorem}[Telg\'arsky]\label{telgarskyth} If TWO has a winning strategy in $\gfin^{\omega}(\mathcal{O},\mathcal{O})$, then the space is $\sigma$-compact.
\end{theorem}
\begin{theorem}[Galvin, Telg\'arsky]\label{galvinth} If TWO has a winning strategy in $\gone^{\omega}(\mathcal{O},\mathcal{O})$, then the space is countable.
\end{theorem}
The $\sone(\mathcal{O},\mathcal{O})$-type of a space $X$, denoted ${\sf tp}_{\sone(\mathcal{O},\mathcal{O})}(X)$ is the least ordinal $\alpha$ such that TWO has a winning strategy in $\gone^{\alpha}(\mathcal{O},\mathcal{O})$, played on $X$. 

These ordinal numbers were introduced in \cite{DG} for the point-open game, and there were called point-open types. Galvin introduced the game $\gone^{\omega}(\mathcal{O},\mathcal{O})$ in \cite{G} and showed that it is the ``dual" of the point-open game. Galvin's techniques can be used to directly translate the results of \cite{DG} to the selection principles context. Taking this  information into account, Theorems \ref{scproduct}, \ref{sigmacompact} and \ref{countable}  can be reformulated as follows, respectively:
\begin{quote} If TWO has a winning strategy in $\gfin^{\omega}(\mathcal{O},\mathcal{O})$ on $X$, and $X$ and $Y$ have $\Sc(\mathcal{O},\mathcal{O})$ then $X\times Y$ has $\Sc(\mathcal{O},\mathcal{O})$.
\end{quote}
\begin{quote} If TWO has a winning strategy in $\gfin^{\omega}(\mathcal{O},\mathcal{O})$ on $X$, and $Y$ has $\sfin(\mathcal{O},\mathcal{O})$ then $X\times Y$ has $\sfin(\mathcal{O},\mathcal{O})$.
\end{quote}
\begin{quote} If TWO has a winning strategy in $\gone^{\omega}(\mathcal{O},\mathcal{O})$ on $X$, and $Y$ has $\sone(\mathcal{O},\mathcal{O})$ then $X\times Y$ has $\sone(\mathcal{O},\mathcal{O})$.
\end{quote}
These reformulations suggest another alternative to generalizing Theorems \ref{scproduct}, \ref{sigmacompact} and \ref{countable}. In the next section we give some results in this direction.

\section{Products when one factor has a strong form of Rothberger's property}

We will use the following equivalences in the proof of Theorem \ref{pot} below:
\begin{enumerate}
\item{$\Sc(\mathcal{O},\mathcal{O}) \Leftrightarrow \Sc(\Omega,\mathcal{O})$ - \cite{lbFilomat}.}
\item{$\sfin(\mathcal{O},\mathcal{O}) \Leftrightarrow \sfin(\Omega,\mathcal{O})$ - \cite{coc1}.}
\item{$\sone(\mathcal{O},\mathcal{O}) \Leftrightarrow \sone(\Omega,\mathcal{O})$ - \cite{coc1}.}
\end{enumerate}

We shall also use equivalent forms of the property $\sfin(\Omega,\mathcal{O}^{gp})$. The \emph{Hurewicz property} of $X$ as originally defined in \cite{wh} by Hurewicz is the statement that for each sequence $(\mathcal{U}_n:n<\infty)$ of open covers of $X$ there is a sequence $(\mathcal{V}_n:n<\infty)$ of finite sets such that for each $n$ $\mathcal{V}_n\subseteq \mathcal{U}_n$, and for each $x\in X$, for all but finitely many $n$ we have $x\in\cup\mathcal{V}_n$. The \emph{Hurewicz game} is the game of length $\omega$ where for each $n$, in the $n$-th inning ONE chooses an open cover $O_n$ of $X$, and TWO responds by choosing a finite set $T_n\subset O_n$. A play $(O_1,\, T_1,\, O_2,\,T_2,\,\cdots,\, O_n,\, T_n,\, \cdots)$ is won by TWO if each $x\in X$ is in all but finitely many of the sets $\cup T_n$. Else, ONE wins the play.

In Theorem 27 of \cite{coc1} it was proved that a space $X$ has the Hurewicz property if, and only if, ONE has no winning strategy in the Hurewicz game (In \cite{coc1} the result is stated for $X$ a set of real numbers, but the argument works in the general situation). Then in Theorem 14 of \cite{coc2} it was proved that in separable metric spaces the Hurewicz property is equivalent to $\sfin(\Omega,\mathcal{O}^{gp})$. Theorem 14 (5) of \cite{coc7} is stated for $\sfin(\Omega,\Lambda^{gp})$, which is equivalent to $\sfin(\Omega,\mathcal{O}^{gp})$.

\begin{lemma}\label{scandhurewicz} 
The following statements are equivalent:
\begin{enumerate}
\item{$X$ has properties $\Sc(\mathcal{O},\mathcal{O})$ and $\sfin(\Omega,\mathcal{O}^{gp})$.}
\item{For each sequence $(\mathcal{U}_n:n<\infty)$ of $\omega$ covers of $X$ there is a sequence $(\mathcal{V}_n:n<\infty)$ such that:
      \begin{enumerate}
      \item{Each $\mathcal{V}_n$ is a finite collection of open sets;}
      \item{Each $\mathcal{V}_n$ is pairwise disjoint;}
      \item{Each $\mathcal{V}_n$ refines $\mathcal{U}_n$;}
      \item{For $m\neq n$, $\mathcal{V}_m\cap\mathcal{V}_n = \emptyset$ and}
      \item{there is a sequence $n_1 < n_2 < \cdots < n_k < \cdots$ of positive integers such that each element of $X$ is in all but finitely many of the sets $\cup(\cup_{n_k\le j< n_{k+1}}\mathcal{V}_n)$.}
      \end{enumerate}
}
\end{enumerate}
\end{lemma}
{\flushleft{\bf Proof:}} To prove (2) $\Rightarrow$ (1) one shows: 
\begin{itemize}
\item{(2) (b), (c) and (e) imply $\Sc(\Omega,\mathcal{O})$, which is equivalent to $\Sc(\mathcal{O},\mathcal{O})$.}
\item{(2) (a), (c), (d) and (e) imply (with some work) $\sfin(\Omega,\mathcal{O}^{gp})$.}
\end{itemize} 
We shall now prove (1) $\Rightarrow$ (2):

Let a sequence $(\mathcal{U}_n:n<\infty)$ of $\omega$-covers of $X$ be given. Define a strategy $\sigma$ of player ONE in the Hurewicz game as follows:

Choose by (1) a sequence $(\mathcal{W}_n:n<\infty)$ such that each $\mathcal{W}_n$ is a refinement of $\mathcal{U}_n$ by a pairwise disjoint family of open sets, and $\bigcup_{n<\infty}\mathcal{W}_n$ is an open cover of $X$. Define:
\[
  \sigma(\emptyset) := \bigcup_{n<\infty}\mathcal{W}_n.
\] 
Let TWO choose a finite set $W_1\subset\sigma(\emptyset)$.

To define $\sigma(W_1)$, define: $n_1 = \min\{n:W_1\subseteq \bigcup_{j\le n}\mathcal{W}_j\}$ and $\epsilon_1 = \min(\{diam(V):V\in W_1\}\cup\{1\})$. Then again by (1) choose for each $n>n_1$ a \underline{new} pairwise disjoint family $\mathcal{W}_n$ of open sets such that $\mathcal{W}_n$ refines $\mathcal{U}_n$, and for each $V\in\mathcal{W}_n$ also $diam(V)<\epsilon_1$, and $\bigcup_{n>n_1}\mathcal{W}_n$ is a cover of $X$. Define:
\[
  \sigma(W_1) = \bigcup_{n>n_1}\mathcal{W}_n.
\]
Let TWO choose a finite set $W_2\subset \sigma(W_1)$.

To define $\sigma(W_1,W_2)$, define: $n_2 = \min\{n>n_1:W_2\subseteq \bigcup_{n_1<j\le n}\mathcal{W}_j\}$ and $\epsilon_2 = \min\{diam(V):V\in W_2\}$. Then again by (1) choose for each $n>n_2$ a \underline{new} pairwise disjoint family $\mathcal{W}_n$ of open sets such that $\mathcal{W}_n$ refines $\mathcal{U}_n$, and for each $V\in\mathcal{W}_n$ also $diam(V)<\epsilon_2$, and $\bigcup_{n>n_2}\mathcal{W}_n$ is a cover of $X$. Define:
\[
  \sigma(W_1,W_2) = \bigcup_{n>n_2}\mathcal{W}_n.
\]
It is now clear how the strategy $\sigma$ of ONE is defined.

By (1) and Theorem 27 of \cite{coc1}, ONE has no winning strategy in the Hurewicz game on $X$. Thus, consider a $\sigma$-play lost by ONE, say
\[
  \sigma(\emptyset),\, W_1,\, \sigma(W_1),\, W_2,\, \sigma(W_1,W_2),\, \cdots,\, W_n, \sigma(W_1,\, W_2,\, \cdots,\, W_n),\, \cdots.
\]
Using the definition of $\sigma$ we find sequences $n_1\, < n_2\, < n_3\, \cdots$ of positive integers and $\epsilon_1 \,> \, \epsilon_2\, > \, \epsilon_3\, > \, \cdots$ of positive real numbers, and families $\mathcal{W}_n$, $n<\infty$ of pairwise disjoint open sets such that:
\begin{enumerate}
\item{$W_1\subseteq\bigcup_{j\le n_1}\mathcal{W}_j$;}
\item{For each $k$, $W_{k+1}\subseteq \bigcup_{n_k<j\le n_{k+1}}\mathcal{W}_j$;}
\item{For each $k$, $\mathcal{W}_k$ refines $\mathcal{U}_k$;}
\item{For each $k$, for each $V\in W_k$, $\epsilon_{k+1}<\epsilon_k\le diam(V) < \epsilon_{k-1}$;}
\item{For each $x\in X$, for all but finitely many $k$, $x\in\bigcup W_k$.}
\end{enumerate}
By (4) we have that for $i\neq j$, $W_i\cap W_j = \emptyset$.

For each $j$ define $\mathcal{V}_j$ as follows:
\begin{itemize}
\item{$j\le n_1$: $\mathcal{V}_j = (W_1\cap\mathcal{W}_j)\setminus(\bigcup_{i<j}\mathcal{V}_i)$;}
\item{$n_k<j\le n_{k+1}$: $\mathcal{V}_j = (W_{k+1}\cap\mathcal{W}_j)\setminus(\bigcup_{i<j}\mathcal{V}_i)$.}
\end{itemize}
The sequence $(\mathcal{V}_n:n<\infty)$ is as required.
$\diamondsuit$

The following lemma is also used in the proof of Theorem \ref{pot} below:
\begin{lemma}\label{omegacoverrefinement} Every $\omega$-cover of $X\times Y$ is refined by one whose elements are of the form $U\times V$ where $U\subset X$ and $V\subset Y$ are open. 
\end{lemma}
{\flushleft{\bf Proof:}} For each finite $F\subset X\times Y$ choose finite sets $F_X$ and $F_Y$ with $F\subseteq F_X\times F_Y$, and an open set $U\in\mathcal{U}_n$ with $F_X\times F_Y\subset U$. Then choose open sets $U_X$ and $U_Y$ with $F_X\subset U_X$ and $F_Y\subset U_Y$, and $U_X\times U_Y\subset U$.  $\diamondsuit$

\begin{theorem}\label{pot} Let $X$ be a metric space with ${\sf tp}_{\sone(\mathcal{O},\mathcal{O})}(X)<\omega^2$.
\begin{enumerate}
\item{If $Y$ has $\Sc(\mathcal{O},\mathcal{O})$ and $\sfin(\Omega,\mathcal{O}^{gp})$, then $X\times Y$ has $\Sc(\mathcal{O},\mathcal{O})$.}
\item{If $Y$ has $\sfin(\Omega,\mathcal{O}^{gp})$, then $X\times Y$ has $\sfin(\mathcal{O},\mathcal{O})$.}
\item{If $Y$ has $\sone(\Omega,\mathcal{O}^{gp})$, then $X\times Y$ has $\sone(\mathcal{O},\mathcal{O})$.}
\end{enumerate}
\end{theorem}
{\flushleft {\bf Proof:}} We prove this by induction on ${\sf tp}_{\sone(\mathcal{O},\mathcal{O})}(X)$. Note that Theorem 5(a) and (b) of \cite{DG} imply that ${\sf tp}_{\sone(\mathcal{O},\mathcal{O})}(X)$ ($<\omega^2$) is a limit ordinal. 

For $\alpha=\omega\cdot 1$ the Galvin-Telg\'arsky theorem implies that $X$ is countable (and thus $\sigma$-compact). By Theorem \ref{scproduct} we have (1), by Theorem \ref{sigmacompact} we have (2) and by Theorem \ref{countable} we also have (3). 

Now suppose we have proven (1), (2) and (3) for all metric spaces $X$ with ${\sf tp}_{\sone(\mathcal{O},\mathcal{O})}(X)\le \omega\cdot n$. Consider a metric space $X$ with ${\sf tp}_{\sone(\mathcal{O},\mathcal{O})}(X) = \omega\cdot (n+1)$. By Theorem 5(b) of \cite{DG} fix a countable set $D\subset X$ such that for each open subset $U\supseteq D$ of $X$, ${\sf tp}_{\sone(\mathcal{O},\mathcal{O})}(X\setminus U)\le \omega\cdot n$.

By the remarks preceding Lemma \ref{scandhurewicz} it suffices to consider only $\omega$-covers of $X\times Y$. Let $(\mathcal{U}_n:n<\infty)$ be a sequence of $\omega$-covers of $X\times Y$. By Lemma \ref{omegacoverrefinement} we may assume that the elements of each $\mathcal{U}_n$ are of the form $U\times V$.  
For each $n$ and for each finite $F\subset D$ put
\[
  \mathcal{S}_{F,n} = \{V:(\exists U)(F\subset U \mbox{ and } U\times V\in\mathcal{U}_n)\}.
\]
For each $V\in\mathcal{S}_{F,n}$ choose a $U_V$ open in $X$ with $F\subset U_V$ and $U_V\times V\in\mathcal{U}_n$. 
Write $D = \cup_{n<\infty} F_n$ where for each $n$ $F_n\subset F_{n+1}$, and $F_n$ is a finite set. 
Now $(\mathcal{S}_{F_n,2^n}:n<\infty)$ is a sequence of $\omega$-covers of $Y$.\\
{\flushleft{\bf Proof of (1):}} Assume that $Y$ has $\Sc(\mathcal{O},\mathcal{O})$ and $\sfin(\Omega,\mathcal{O}^{gp})$:\\
Choose by Lemma \ref{scandhurewicz} for each $n$ a finite pairwise disjoint set $\mathcal{V}_{F_n,2^n}$ of open sets refining $\mathcal{S}_{F_n,2^n}$ such that for $n\neq m$ we have $\mathcal{V}_{F_n,2^n}\cap \mathcal{V}_{F_m,2^m}=\emptyset$. Also choose a sequence $n_1< n_2 < \cdots < n_k<\cdots$ of positive integers such that each element of $Y$ is in all but finitely many of the sets $\bigcup(\cup_{n_k\le j<n_{k+1}}\mathcal{V}_{F_j,2^j})$.

For each $k$, put $U_k = \cap\{U_V:\, V\in \cup_{n_k\le j<n_{k+1}}\mathcal{V}_{F_j,2^j}\}$. As finite intersection of open sets, $U_k$ is open. Also, $U_k\supseteq F_{n_k}$. Thus, as for each $k$ $F_{n_k}\subseteq F_{n_{k+1}}$ we have: $\{U_k:k<\infty\}$ is a $\gamma$-cover of $D$.

Each $G_m = \cup_{j>m}U_j$ is an open subset of $X$ and contains $D$. Each $X_m = X\setminus G_m$ is a closed set in $X$ and has ${\sf tp}_{\sone(\mathcal{O},\mathcal{O})}(X_m) \le \omega\cdot n$. Put $G = \cap_{m<\infty}G_m$. By Theorem 4 (e) of \cite{DG}, ${\sf tp}_{\sone(\mathcal{O},\mathcal{O})}(X\setminus G) = {\sf tp}_{\sone(\mathcal{O},\mathcal{O})}(\bigcup_{m<\infty}X_m) \le \omega\cdot n$. Moreover, $\{U_n:n<\infty\}$ is a large cover of $G$. Thus for each $j$ in $[n_k,n_{k+1})$ the set $\{U_k\times V:V\in\mathcal{V}_{F_j,2^j}\}$ is a pairwise disjoint refinement of $\mathcal{U}_{2^j}$,  and $\bigcup_{k<\infty}(\cup_{n_k\le j<n_{k+1}} \{U_k\times V:V\in \mathcal{V}_{F_j,2^j}\})$ is a large cover of $G \times Y$. 

And as ${\sf tp}_{\sone(\mathcal{O},\mathcal{O})}(X\setminus G) \le \omega\cdot n$, the induction hypothesis implies that there is a sequence $({\mathcal V}_j:j<\infty \mbox{ and }(\forall n)(j\neq 2^n))$ such that for each such $j$, $\mathcal{V}_j$ is a pairwise disjoint family which refines $\mathcal{U}_j$ and the union of these $\mathcal{V}_j$ is an open cover of $(X\setminus G)\times Y$. 

Thus we obtained a sequence of sets $\mathcal{V}_n$, each a pairwise disjoint family refining $\mathcal{U}_n$, such that $\cup_{n<\infty}\mathcal{V}_n$ is an open cover of $X\times Y$.

{\flushleft{\bf Proof of (2):}} Assume that $Y$ has the property $\sfin(\Omega,\mathcal{O}^{gp})$:\\
Choose for each $n$ a finite subset $\mathcal{V}_{F_n,2^n}$ of $\mathcal{S}_{F_n,2^n}$ such that each element of $Y$ is in all but finitely many of the sets $\cup\mathcal{V}_{F_n,2^n}$. For each $n$, put $U_{2^n} = \cap\{U_V:\, V\in \mathcal{V}_{F_n,2^n}\}$. As finite intersection of open sets, $U_{2^n}$ is open. Also, $\{U_{2^n}:n<\infty\}$ is a $\gamma$-cover of $D$. For each $n$ put $G_n = \cup_{j>n}U_{2^j}$, an open subset of $X$, containing $D$. Then each $X_m = X\setminus G_m$ is a closed set in $X$ and has ${\sf tp}_{\sone(\mathcal{O},\mathcal{O})}(X_m) \le \omega\cdot n$. By Theorem 4 (e) of \cite{DG}, ${\sf tp}_{\sone(\mathcal{O},\mathcal{O})}(X\setminus(\cap_{m<\infty}G_m)) = {\sf tp}_{\sone(\mathcal{O},\mathcal{O})}(\bigcup_{m<\infty}X_m) \le \omega\cdot n$. Moreover, $\{U_{2^n}:n<\infty\}$ is a large cover of $G = \cap_{m<\infty}G_m$. Thus $\cup_{n<\infty}\{U_{2^n}\times V:V\in\mathcal{V}_{F_n,2^n}\}$ is a large cover of $G \times Y$. And as ${\sf tp}_{\sone(\mathcal{O},\mathcal{O})}(X\setminus G) \le \omega\cdot n$, the induction hypothesis implies that there is a sequence $({\mathcal V}_j:j<\infty \mbox{ and }(\forall n)(j\neq 2^n))$ of finite sets such that for each such $j$, $\mathcal{V}_j\subseteq \mathcal{U}_j$, and the union of these $\mathcal{V}_j$ covers $(X\setminus G)\times Y$. Thus we obtain a sequence of finite sets $\mathcal{V}_n\subset \mathcal{U}_n$ such that $\cup_{n<\infty}\mathcal{V}_n$ covers $X\times Y$.\\

{\flushleft{\bf Proof of (3):}} Assume $Y$ has property $\sone(\Omega,\mathcal{O}^{gp})$:\\
Choose for each $n$ a set $V_{F_n,2^n} \in \mathcal{S}_{F_n,2^n}$ such that $\mathcal{V} = \{V_{F_n,2^n}:n<\infty\}$ is a  groupable cover of $Y$. Now we can chooses sequences $\ell_1 < m_1 < \ell_2 < m_2 < \cdots < \ell_k < m_k <\cdots$ such that for each $y\in Y$, for all but finitely many $k$ we have $y \in \cup_{\ell_k \le j\le m_k}V_{F_j,2^j}$. For each $k$ put $W_k = \cap_{\ell_k \le j\le m_k}U_{F_j,2^j}$, a finite intersection of open sets containing $F_{\ell_k}$. Then $\{W_k:k<\infty\}$ is a $\gamma$-cover of $D$, and as above we find a ${\sf G}_{\delta}$ set $G\subset X$ such that each element of $G$ is in infinitely many $W_k$, and so $G\times Y$  is covered by the sets $U_{F_j,2^j}\times V_{F_j,2^j}$ selected from $\mathcal{U}_{2^j}$ for all these $j$.  And since ${\sf tp}_{\sone(\mathcal{O},\mathcal{O})}(X\setminus G) \le \omega\cdot n$ apply the induction hypothesis to $(X\setminus G)\times Y$.
$\diamondsuit$

A set of real numbers is a \emph{Lusin} set if it is uncountable but its intersection with each nowhere dense set is countable. And  a set of reals is a \emph{Sierpi\'nski} set if it is uncountable but its intersection with each Lebesgue measure zero set is countable. The Continuum Hypothesis implies that these exist. Rothberger proved that conversely, if there are both Sierpi\'nski and Lusin sets of cardinality $2^{\aleph_0}$, then $2^{\aleph_0} = \aleph_1$ - \cite{rch}. Lusin sets have $\sone(\mathcal{O},\mathcal{O})$ but their finite powers can fail to have $\sfin(\mathcal{O},\mathcal{O})$. Sierpi\'nski sets have $\sone(\Gamma,\Gamma)$, and thus $\sfin(\Omega,\mathcal{O}^{gp})$, but their finite powers may fail to have $\sfin(\mathcal{O},\mathcal{O})$. 

\begin{corollary}\label{LusinSierpinski} For $X$ a Lusin set and $Y$ a Sierpi\'nski set, $X\times Y$ has $\sfin(\mathcal{O},\mathcal{O})$.
\end{corollary}

\begin{corollary}\label{gammaLusin} For $X$ a Lusin set and $Y$ a $\gamma$-set, $X\times Y$ has $\sone(\mathcal{O},\mathcal{O})$.
\end{corollary}

\section{Remarks}

The bound $\omega^2$ in the hypothesis of Theorem \ref{pot} should probably be $\omega_1$. In fact, the proof of this theorem applies to all the sets of reals constructed  in \cite{Baldwin} with the aid of the Continuum Hypothesis: Each countable limit ordinal is ${\sf tp}_{\sone(\mathcal{O},\mathcal{O})}(X)$  for some such $X$.
\begin{conjecture}\label{bound} For every metric space $X$ with  ${\sf tp}_{\sone(\mathcal{O},\mathcal{O})}(X)$ countable:
\begin{enumerate}
\item{If $Y$ has ${\Sc}(\mathcal{O},\mathcal{O})$ and $\sfin(\Omega,\mathcal{O}^{gp})$, then $X\times Y$ has ${\Sc}(\mathcal{O},\mathcal{O})$.}
\item{If $Y$ has property $\sfin(\Omega,\mathcal{O}^{gp})$, then $X\times Y$ has $\sfin(\mathcal{O},\mathcal{O})$.}
\item{If $Y$ has property $\sone(\Omega,\mathcal{O}^{gp})$, then $X\times Y$ has $\sone(\mathcal{O},\mathcal{O})$.}
\end{enumerate}
\end{conjecture}
There are models of Set Theory in which ${\sf tp}_{\sone(\mathcal{O},\mathcal{O})}(X)$ is $\omega$ or $\omega_1$ for infinite sets $X$ of real numbers: Any model in which the {\tt Borel Conjecture} or {\tt Martin's Axiom plus negation of the Continuum Hypothesis} holds will do. The result regarding Martin's Axiom is due to Fremlin - see \cite{DG}. Conjecture \ref{bound} is for trivial reasons true in these models.  

Since a product of Lusin sets need not be $\sfin(\mathcal{O},\mathcal{O})$ (\cite{coc2} Theorems 2.6 and 3.1) and Lusin sets are concentrated on countable subsets of themselves, it follows that we cannot replace the Hurewicz property $\sfin(\Omega,\mathcal{O}^{gp})$ in Theorem \ref{pot} (2) with merely the Menger property $\sfin(\mathcal{O},\mathcal{O})$, nor with the property of being concentrated on a countable (or $\sigma$-compact) subset of itself. 

Can Theorem \ref{pot} (2) be proved by merely assuming that TWO has a winning strategy in $\gfin^{\alpha}(\mathcal{O},\mathcal{O})$ for some countable ordinal $\alpha$? The answer is ``no": It was pointed out after Theorem 2.11 of \cite{coc2} that the Continuum Hypothesis implies that there is a Sierpi\'nski set $S$ such that $S\times S$ does not have the property $\sfin(\mathcal{O},\mathcal{O})$. But Sierpi\'nski sets have $\sone(\Gamma,\Gamma)$ (\cite{coc2}, Theorem 2.10), and so $\sfin(\Omega,\mathcal{O}^{gp})$. Moreover:
\begin{proposition}\label{MengerSierpinski} If $S\subset\reals$ is a Sierpinski set, then TWO has a winning strategy in $\gfin^{\omega+\omega}(\mathcal{O},\mathcal{O})$.
\end{proposition} 
{\flushleft{\bf Proof:}} TWO's strategy is to make sure that in the $n$-th inning the outer measure of $([-n,n]\cap S)\setminus(\bigcup T_n)$ is less than $\frac{1}{2^n}$. After $\omega$ innings the measure of the piece of $S$ not yet covered by TWO is $0$, and thus countable. In the remaining $\omega$ innings TWO covers this countable set. $\diamondsuit$ 

Can Theorem \ref{pot} (1) be proved by merely assuming that TWO has a winning strategy in $\gfin^{\alpha}(\mathcal{O},\mathcal{O})$ for some countable ordinal $\alpha$? We don't know. The following problem presents a good test case:
\begin{problem}
If $X$ is a Sierpi\'nski set and $Y$ has ${\Sc}(\mathcal{O},\mathcal{O})$ and $\sfin(\Omega,\mathcal{O}^{gp})$, then does $X\times Y$ have property ${\Sc}(\mathcal{O},\mathcal{O})$?

\end{problem}

\bigskip

\begin{center}
  Addresses
\end{center}

\medskip

\begin{flushleft}
Liljana Babinkostova                      \\
Boise State University                    \\
Department of Mathematics                 \\
Boise, ID 83725 USA                       \\
e-mail: liljanab@math.boisestate.edu      \\
\end{flushleft}

\begin{flushleft}
Marion Scheepers                    \\
Department of Mathematics           \\
Boise State University              \\
Boise, Idaho 83725 USA              \\
e-mail: marion@math.boisestate.edu   \\
\end{flushleft}

\end{document}